\documentclass[12pt]{amsart}
\usepackage{amsmath}
\usepackage{latexsym}
\usepackage{amsthm}
\usepackage{amscd}
\usepackage{amsfonts}
\usepackage{graphicx}
\newcommand{\R}{\mathbb R}

\newcommand{\fg}{\mathfrak{g}}

\newcommand{\fa}{\mathfrak{a}}
\newcommand{\SL}{\mathop{\mathrm{SL}}}

\newcommand{\GL}{\mathop{\mathrm{GL}}}

\newcommand{\norm}[1]{\left\Vert#1\right\Vert}

\newcommand{\abs}[1]{\left\vert#1\right\vert}
\newcommand{\set}[1]{\left\{#1\right\}}

\newcommand{\rip}[2]{\left(#1\, ,\, #2 \right)}

\newcommand{\ad}{{\operatorname{ad}}}

\theoremstyle{plain}
\newtheorem{thm}{Theorem}[section]

\newtheorem{lem}[thm]{Lemma}
\newtheorem{prop}[thm]{Proposition}

\theoremstyle{definition}

\newtheorem{rem}[thm]{Remark}

\def\tr{\operatorname{tr}}
\def\Tr{\operatorname{Tr}}
\def\det{\operatorname{det}}
\def\Det{\operatorname{Det}}

\numberwithin{equation}{section}

\begin{document}
\title[Recursion Relations for Laguerre Functions ]
{Differential Recursion Relations for Laguerre Functions on
Symmetric Cones}
\author{Michael Aristidou \and Mark Davidson\and Gestur \'Olafsson }
\address{Department of Mathematics, Louisiana State University, Baton Rouge, LA\ 70803, USA}
\email{aristido@math.lsu.edu} \email{davidson@math.lsu.edu}
\email{olafsson@math.lsu.edu}
\thanks{M. Aristidou was partially supported by a departmental Focussed Research
Award and by an NSF grant DMS-0139783. G. \'Olafsson was supported
by NSF grants DMS-0139783 and DMS-0402068.}

\dedicatory{Dedicated to Jacques Faraut } 


\begin{abstract}
Let $\Omega$ be a symmetric cone and $V$ the corresponding simple
Euclidean Jordan algebra. In \cite{ado,do,do04,doz2} we considered
the family of generalized Laguerre functions on $\Omega$ that
generalize the classical Laguerre functions on $\mathbb{R}^+$. This
family  forms an orthogonal basis for the subspace of $L$-invariant
functions   in $L^2(\Omega,d\mu_\nu)$, where $d\mu_\nu$ is a certain
measure on the cone and where $L$ is the group  of linear
transformations on $V$ that leave the cone $\Omega$ invariant and
fix the identity in $\Omega$. The space $L^2(\Omega,d\mu_\nu)$
supports a highest weight representation of the group $G$ of
holomorphic diffeomorphisms that act on the tube domain
$T(\Omega)=\Omega + iV.$ In this article we give an explicit formula
for the action of the Lie algebra of $G$ and via this action
determine  second order differential operators which give
differential  recursion relations for the generalized Laguerre
functions generalizing the classical creation, preservation, and
annihilation relations for the Laguerre functions on $\mathbb{R}^+$.
\end{abstract}
\maketitle

\section*{Introduction}\label{s0}
\noindent It is a general understanding that special functions   are
closely related to representation theory of special Lie groups.
Special functions are realized as coefficient functions of the
representation and the action of the Lie algebra is used to derive
differential equations and recursion relations satisfied by those
functions. Standard references to this philosophy are the works of
Vilenkin and Klimyk  \cite{v68,vk91}. We also refer the interested
reader to the text \cite{d80}, the recent text \cite{aar}, and the
work of T. Koornwinder. The present article reflects these general
philosophies. In particular,  we conclude our work on the connection
between generalized Laguerre functions, highest weight
representations and Jordan algebras, \cite{ado,do,do04,doz1,doz2}.
The classical Laguerre functions $\ell_n^\lambda$ form an orthogonal
basis for the Hilbert space $L^2(\mathbb{R}^+,x^{\lambda -1}dx)$,
$\lambda >0$. As far as we have been able to trace, the first
generalizations of the Laguerre functions and polynomials is from
1935  in the work of F. Tricomi \cite{t35}. Later, in 1955,  C. S.
Herz \cite{h55} considered generalized Laguerre functions in the
context of Bessel functions on the space of complex $m\times
m$-matrices. The Laguerre polynomials are defined on the cone of
positive definite complex matrices   in terms of the generalized
hypergeometric functions, also introduced in the same article.
Other realizations of the Laguerre functions, using the Laplace
transform, were also derived. The motivation was to construct a
complete set of eigenfunctions for the Hankel transform. The
generalization to all symmetric cones using Euclidean Jordan
algebras was achieved almost 40 years later in the beautiful book by
J. Faraut and A. Koranyi \cite{fk}. Here the Laguerre polynomials
were defined in terms of certain polynomials $\psi_{\mathbf{m}}(x)$
invariant under the action of a maximal compact subgroup $L$ leaving
the cone invariant and fixing the identity $e$:
$$L^\nu_\mathbf{m}(x)=(\nu )_{\mathbf{m}}\sum_{|\mathbf{n}|\le \mathbf{m}}
 \binom{{\mathbf
m}}{{\mathbf n}} \frac{1}{(\nu)_{\mathbf n} }\psi_{\mathbf n}(-x)\, ,$$
c.f. Section \ref{s4}. The Laguerre functions are defined as
$$\ell_\mathbf{m}^\nu (x)=e^{-\tr x}L^\nu_\mathbf{m}(2x ),$$
where $\tr$ is the trace in the corresponding Jordan algebra. It was
shown that the Laguerre functions were orthogonal and eigenfunctions
of the Hankel transform. Later, F. Ricci and A. Tabacco
constructed a system of differential operators, in the context of
the Jordan algebra of Hermitian symmetric matrices and real
symmetric matrices, having the Laguerre functions as eigenfunctions
with distinct eigenvalues, c.f. \cite{RiV}. In the simplest case
this differential operator is nothing but the Laguerre differential
operator. None of these works, however, considers the generalized
differential recursion relations that correspond to raising and
lower operators satisfied by the Laguerre functions.

The first time that the Laguerre polynomials were directly related
to representation theory was in \cite{v68} where they were shown to
be coefficient functions of representations of the group
$$\left\{\left(
\begin{matrix}1 & a & b\cr 0 & c & d\cr 0 & 0 &1\cr
\end{matrix}
\right)\mid a,b,c,d\in \mathbb{R}, c\not= 0\right\}\, . $$ Later B.
Kostant and N. Wallach used the recursion relations that exist
amongst the Laguerre functions to construct a highest weight
representation and subsequently study Whittaker vectors for some
special representations \cite{k00,w00}. In \cite{k00} the
differential equations and recursion relations for the Laguerre
functions were used to give a realization of the highest weight
representations of $\widetilde{\mathrm{SL}(2,\mathbb{R})}$, the
universal covering group of $\mathrm{SL}(2,\mathbb{R})$. The
opposite point of view was taken in \cite{doz1} where the authors
showed how one can derive those classical relations using a highest
weight representation  and the Laplace transform. The classical
relations were given as the action of special elements in the Lie
algebra acting as second order differential operators on functions
 on $\mathbb{R}^+$.

The connection to the construction in \cite{fk} was established in
\cite{doz2} where the generalized Laguerre functions were shown to
be not only orthogonal but also a basis of the  space of
$L$-invariant functions in the highest weight module realized in
$L^2(\Omega ,d\mu_\nu)$, where $\Omega$ is a symmetric cone, and
$d\mu_\nu$ is a certain quasi-invariant measure on $\Omega$. Using a
certain $L$-invariant element in the Lie algebra, the authors showed
that the Laguerre functions satisfy a first order differential
recursion relation involving the Euler operator (c.f. Theorem 7.9 in
\cite{doz2}).  The terms in this relation involve a raising and
lowering of indices that parameterize the Laguerre functions. In the
context of a highest weight representation one deduces that the
Euler operator is made up of a creation and annihilation operator
derived from the action of the Lie algebra. However, no attempt was
made to derive an explicit form of these operators until we
considered the special cases of the cones of Hermitian symmetric
matrices and real symmetric matrices in \cite{do,ado}, respectively.
In this article we generalize those results to the Laguerre
functions related to all symmetric cones. The tools are again
highest weight representations and Jordan algebras. The main results
are the explicit formulas for the action of the Lie algebra in the
realization of the highest weight space $L^2(\Omega,d\mu_\nu)$  and
then the restriction to the subalgebra of $L$-invariants which
results in the differential equations and recursion relations in
terms of explicitly constructed differential operators.

If $\mathfrak{g}$ is simple, then the subalgebra of $L$-invariants
in $\mathfrak{g}_\mathbb{C}$ is isomorphic to
$\mathfrak{sl}(2,\mathbb{C})$. It should be noted that such a three
dimensional Lie algebra  of differential operators has shown up in
several places in the literature. We would like to mention its
important role in the study of the Huygens' principle
\cite{bo05,h80,ht92}, in representation theory \cite{km05} (and the
references therein), and in the  theory of special functions
\cite{r98}.

One cannot downplay the essential role that Jordan algebras play in
establishing and expressing many of the fundamental results obtained
about orthogonal families of special functions defined on symmetric
cones. Nevertheless, the theory of highest weight representations
adds fundamental new results not otherwise easily obtained. In
short, our philosophy is  that there is a strong interplay between
Jordan algebras, highest weight representations, and special
functions which has not been fully exploited.

The starting point in this project has been the representation
theory, wherein the Laguerre polynomials form an orthogonal family
of functions invariant under a group action. However, the Laguerre
polynomials have also been introduced in the literature using
several variable Jack polynomials \cite{bf97,d97,l91}. We would like
to thank M. R\"osler for pointing these references out to us. To
explain the connection, a little more notation is needed. Let $J$ be
an irreducible Euclidean Jordan algebra of rank $r$. Let $c_1,\ldots
,c_r\in J$ be a Jordan frame, $\mathfrak{a}=\bigoplus_{j=1}^r\R c_j$
and $e=c_1+\ldots + c_r$. Let $\Omega=\{x^2\mid x\in J\, \, \text{
and $x$ regular}\}$ be the standard symmetric cone in $J$. Let
$H=\{g\in \GL (J)\mid g\Omega =\Omega\}_o$ and $L$ the maximal
compact subgroup of $H$ fixing $e$. Then the Laguerre functions and
polynomials are $L$-invariant functions on $\Omega$. Let
$$\Omega_1=\fa \cap \Omega\simeq (\R^+)^r\, .$$
Then $\Omega =L\cdot \Omega_1$ and therefore the Laguerre
polynomials and functions are uniquely determined by their
restriction to $\Omega_1$. Thus,  the Laguerre polynomials can also
be defined as polynomials on $\Omega_1$ or the vector space $\fa$,
invariant under the Weyl group $W_H=N_L(\fa )/Z_L (\fa )$. This is
the way the Laguerre polynomials are defined in the above
references.

In the case of symmetric matrices, this
boils down to the fact that each symmetric matrix can be
diagonalized. Thus
$$\Omega_1=\{d(\lambda_1,\ldots ,\lambda_n)\mid \lambda_j>0\}\, $$
and the Laguerre polynomials can be viewed as polynomials in the
eigenvalues, invariant under permutations.

The article is organized as follows. The necessary tools from Jordan
algebra theory are introduced in Section \ref{s1}. In Section
\ref{s2} we introduce the tube domain $T(\Omega)=V+i\Omega$, where
$V$ is a simple Euclidean Jordan algebra. The main part of this
section is devoted to describing the Lie algebra of
$G=G(T(\Omega))_o$, where $G(T(\Omega ))$ is the group of
holomorphic automorphisms of $T(\Omega )$ and the subscript ${}_o$
stands for the connected component of the identity. The final result
is the description of the $L$-invariant elements in $\mathfrak{g}$.
Most of this material can be found in \cite{fk}.

In Section \ref{s3} we introduce the highest weight representations
and give an explicit realization of those representations in
$L^2(\Omega ,d\mu_\nu )$ using the second order Bessel differential
operator introduced in \cite{fk}:
\begin{equation*}
\mathcal{B }_\nu=P\left(\frac{\partial}{\partial x}\right)x + \nu
\frac{\partial}{\partial x},
\end{equation*}
where $P(a)$ denotes the quadratic representation of $V$.
For $w\in V_\mathbb{C}$ we define the differential operator
$\mathcal{B}_{\nu,w}$ by
\begin{equation*}
\mathcal{B}_{\nu ,w}f(x)=\rip{\mathcal{B}_\nu f(x)}{w}\, .
\end{equation*}
Then the following holds:
\medskip

\noindent \textbf{Theorem \ref{th-3.4}.} \textit{ Let $f\in
L^2(\Omega,d\mu_\nu)^\infty$. The representation $\lambda_\nu$ of
${\mathfrak g}$ is described as follows:}\smallskip

$\begin{array}{lll} (1) &\lambda_\nu(X(iu,0,0))f(x)= \tr(iux)f(x),&
X(iu,0,0)\in
{\mathfrak n}^+,\\
(2)&\lambda_\nu(X(0,T,0))f(x)= \frac{\nu}{p}\Tr(T)f(x)+
D_{T^tx}f(x),&
 X(0,T,0)\in {\mathfrak h},\\
(3)&\lambda_\nu(X(0,0,iv))f(x) = -\mathcal{B}_{\nu,iv} f(x),&
X(0,0,iv)\in {\mathfrak  n}^-.
\end{array}$

\medskip \noindent
Here $\mathfrak{n}^+, \mathfrak{h}$, and $\mathfrak{n}^-$ are
certain subalgebras of $\mathfrak{g}$ such that
$\mathfrak{g}=\mathfrak{n}^+\oplus
\mathfrak{h}\oplus\mathfrak{n}^-$.

In Section \ref{s4} we introduce the Laguerre functions and finally,
in Section \ref{s5} we use Theorem \ref{th-3.4} to derive explicit
second order differential operators such that one of them has the
Laguerre functions as eigenfunctions and the two others are,
respectively, an annihilator operator and a creation operator.
Setting $B_\nu = \mathcal{B}_{\nu ,e}$ we have:
\medskip

\noindent
\textbf{Theorem \ref{th-5.2}.}
\textit{The Laguerre functions are related by the following differential
recursion relations:
\begin{enumerate}
\item $(-\tr x + B_\nu )\ell_{\mathbf m}^\nu
  (x)=-(r\nu +2 \abs{{\mathbf m}})\ell_{\mathbf m}^\nu (x)$,
\item $(\tr x +r\nu +2D_x +B_\nu )\ell_{\mathbf m}^\nu(x)=-2\sum_{j=1
}^{r}\binom{\mathbf m}{{\mathbf
m}-\gamma_j}(m_j-1+\nu-(j-1)\frac{d}{2})\ell_{{\mathbf
m}-\gamma_j}^\nu(x)$,
\item $(\tr x -r\nu -2D_x +B_\nu )\ell_{\mathbf
m}^\nu(x)=-2\sum_{j=1}^{r }c_{\mathbf m}(j)\ell_{{\mathbf
m}+\gamma_j}^\nu(x)$,
\end{enumerate}
where the constants $c_{\mathbf m}(j)$  are defined by
$$c_{\mathbf m}(j)=\Pi_{k\ne
j}\frac{m_k-m_j-\frac{d}{2}(k-j+1)}{m_k-m_j-\frac{d}{2}(k-j)}\, .$$
}

\section{Jordan Algebras and Symmetric Cones}\label{s1}
\noindent
In this section we will set down the notation
and basic results concerning Jordan
algebras and symmetric cones used for the remainder of this paper.
  We have tried to keep the notation consistent with the
text by Faraut and Koranyi (c.f. \cite{fk}).    For proofs of
results mentioned below see this text.

Let $V$ be a real   Jordan algebra.   This means that $V$ is a real
vector space with  a bilinear commutative product $(a,b)\mapsto a b$
such that $a^2  (a b)=a(a^2  b)$. In general, a Jordan algebra is
not associative.  Let  $L(a)$ denote left multiplication  by $a$ on
$V$. Thus $L(a)x=ax$. Since the product is bilinear, $L(a)$ is a linear operator on
$V$. The multiplicative property given above is equivalent to
$[L(a),L(a^2)]=0$, for all $a\in V$.

Let
$x\in V$ and let ${\mathbb R}[x]$ be the algebra generated by $x$.
The \textit{rank} of $V$, $r$, is defined by
$$r= \text{max}\set{\dim \mathbb R[x] \mid  x\in V}\, .$$
An element $x\in V$ is
\textit{regular} if $\dim \mathbb R[x]=r$.
The set of regular elements is open and dense in $V$.
Suppose $x$ is regular. We define $\tr(x)$ and $\det(x)$ as follows:
\begin{eqnarray*}
\tr(x) &=& \Tr(L(x)\vert_{{\mathbb R}[x]})\\
\det(x)&=& \Det(L(x)\vert_{{\mathbb R}[x]}),
\end{eqnarray*} where $\Tr$ and $\Det$ are the usual trace and
determinant of a linear operator.  It is not hard to show that
$\tr(x)$ and $\det(x)$ are polynomial functions in $x$ and hence
have polynomial extensions to all of $V$ and $V_{\mathbb{C}}$.

Throughout, we will assume that $V$ is finite dimensional with
dimension $n$ and contains a multiplicative identity $e$.   Let
$V_\mathbb{C}=V\otimes \mathbb{C}$ be the complexification of $V$.
An element $x\in V_{\mathbb{C}}$ is said to be \textit{invertible}
if there is a $y\in {\mathbb{C}[x]} $ such that $x y=e$. The element
$y$ is necessarily unique, it is called the \textit{inverse} of $x$,
and denoted $x^{-1}$. We let $\textsc{i}$ denote the inversion map
on the set of invertible elements: $\textsc{i}(x)=x^{-1}$.

The \textit{quadratic representation} $P$ of $V$ is defined by
$$P(a)=2L(a)^2 - L(a^2)$$
and plays a pivotal role in all that follows. If $F:V\to V$ is
a differentiable map, we denote by $DF:V\to \mathrm{End}(V)$
the derivative of $F$. For $u,x\in V$ we set
$$D_uF(x)=DF(x)u=\lim_{t\to 0}\frac{F(x+tu)-F(x)}{t}\, .$$
\begin{lem}\label{le-DI}
An element $x$ is invertible if and only if $P(x)$ is invertible as
a linear operator on $V$. In this case
\begin{eqnarray*}
P(x)(x^{-1})&=&x\\
P(x)^{-1}&=&P(x^{-1}).
\end{eqnarray*}
The set of invertible elements is an open set in $V$ given by
$\set{x\in V \vert \Det(P(x))\ne 0}$. The derivative of the
inversion map, $\textsc{i}$, is given by $$D
\textsc{i}(x)=-P(x)^{-1},$$ and, in particular, for all $u\in V$, we
have
$$D_u\textsc{i}(x)=-P(x)^{-1}u.$$
\end{lem}

The polarization of $P$ is given by
\begin{eqnarray*}
P(x,y)&:=&\frac{1}{2}D_y(P(x))\\
&=& \frac{1 }{2}(P(x+y)-P(x)-P(y))\\
&=& L(x)L(y)+L(y)L(x)-L(xy).
\end{eqnarray*}

A real Jordan algebra is \textit{Euclidean} if there is an associative inner
product on $V$. In other words, there is an inner product
$\rip{\cdot}{\cdot\,}$ satisfying
$$\rip{x u}{v}=\rip{u}{x v}\, ,$$
for all $x, u, v \in V$. This is equivalent to saying that $L(x)$ is
symmetric for all $x\in V$. A Jordan algebra is simple if there are
no nontrivial ideals.

\begin{prop}
Suppose $V$ is a simple Euclidean Jordan algebra of dimension $n$
and rank $r$. For $x,y\in V$ we
have
\begin{eqnarray*}
\Tr L(x)&=& \frac{n}{r}\tr(x),\\
\Det P(x)&=&(\det x) ^{2n/r},\\
\det(P(y)x)&=& (\det y)^2 \det x.
\end{eqnarray*}
\end{prop}

Henceforth, we will assume
$V$ is a simple Euclidean Jordan algebra of dimension $n$ and rank
$r$. Let $\Omega$ be the interior of the set of all squares $x^2,
x\in V$. Let $G(\Omega)$ be the group of all invertible linear
transformations on $V$ which leave $\Omega$ invariant. We will
also use the notation $H=G(\Omega )_o$, where the subscript ${}_o$ denotes
the connected component containing the identity.

\begin{prop}
The set $\Omega$ is a symmetric cone.  This means that  $\Omega$ is
an open convex cone in $V$, self-dual in the sense that
$$\Omega=\set{y\in V \vert \rip{x}{y}>0, \forall x \in \overline{\Omega}
\;\backslash\set{\,0}}, $$
and  $G(\Omega)$ and $H$
acts transitively on $\Omega$.  Furthermore, $\Omega$ is the
connected component of $e$ in the set of invertible elements of $V$
and
$$\Omega=\set{x\in V \vert L(x) \text{ is positive definite}}.$$
\end{prop}

\section{The Tube Domain $T(\Omega)$}\label{s2}
\noindent Let $V$ be a simple Euclidean Jordan algebra and
$T(\Omega)=\Omega + iV$.\footnote{We choose the right half plane for
$T(\Omega)$ instead of the upper half plane, $V+i\Omega$,  given in
\cite{fk}, for example, and usually referred to as the Siegel upper
half plane.} We note that $T(\Omega )$ is contained in the set of
invertible elements in $V_{\mathbb{C}}$ and $\textsc{i}: z\mapsto
z^{-1}$ is an involutive holomorphic automorphism of $T(\Omega )$
having $e$ as its unique fixed point, c.f. \cite{fk}, Theorem X.1.1.
We note that $V_{\mathbb C}$ is a complex Jordan algebra. The
multiplication, trace, and determinant formulas all extend from $V$
to $V_{\mathbb C}$ in the usual way.  We extend the bilinear form
$\rip{\cdot}{\cdot}$ on $V$ to a complex bilinear form on
$V_{\mathbb C}$ and denote it in the same way.

Let $G(T(\Omega))$ be the group of holomorphic automorphisms of
$T(\Omega)$ and $G=G(T (\Omega ))_o$. We describe elements in $G$ as
follows: Let $(iu,T,iv)\in iV \times H\times iV$ and define
\begin{eqnarray*}
\tau_{iu}(z)&=&z+iu\\
\rho_T(z)&=&Tz\\
\sigma_{iv}(z)&=&(z^{-1}+iv)^{-1}\, .
\end{eqnarray*}
 We observe that
$$\sigma_{iv}=\textsc{i}\tau_{iv}\textsc{i}^{-1}$$
Let $N^+=\set{\tau_{iu}\mid u\in V}$ and $N^-=\set{\sigma_{iv}\mid
v\in V}$.  We identify $H$ with $\set{\rho_T\mid T\in H}$.  It is
well known that the map, $(iu,T,iv)\mapsto \tau_{iu} \rho_T
\sigma_{iv}$, is a diffeomorphism of $N^+\times H\times  N^-$ onto
an open dense subset of $G$. Furthermore, if $K=G_e$, the stabilizer
of $e$, i.e., the set of all $g\in G$ such that $ge=e$, then
$$G=N^+HK\, ,$$ c.f. \cite[p.
205-207]{fk} for details. We set $L=K\cap H$ and note that $K$ is a
maximal compact subgroup of $G$ and $L$ is a maximal compact
subgroup of $H$.

Let ${\mathfrak n^+}$, ${\mathfrak n^-}$, ${\mathfrak h}$,
${\mathfrak g}$, and ${\mathfrak k}$ be the Lie algebras
corresponding to $N^+$, $N^-$, $H$, $G$, and $K$, respectively. The
one parameter subgroups

$$\begin{array}{lll}
z \mapsto z+itu & \in N^+,  & (u\in V); \\[.5ex]
z \mapsto \exp(tT)z &\in H, & (T\in {\mathfrak h});\\[.5ex]
z \mapsto  (z^{-1} +itv)^{-1} & \in N^-, & (v\in V).\\[1ex]
\end{array}$$ induce the corresponding vector fields
$$\begin{array}{ll}
X(z)= iu &\in {\mathfrak n^+},\\[.5ex]
X(z)= Tz & \in {\mathfrak h},\\[.5ex]
X(z)=-P(z)iv & \in {\mathfrak n^-}.\\[1ex]
\end{array}$$
As $\mathfrak{g}=\mathfrak{n}^+\oplus \mathfrak{h}\oplus
\mathfrak{n}^-$ it follows, that every vector field ${\mathfrak g}$
is of the form $$X(z)=iu + T(z) -P(z)(iv)$$ and we will denote it by
the triple $X(iu,T,iv)$. For $x,y\in V_\mathbb{C}$ set
$$x\Box y=L(xy)+[L(x),L(y)]\, .$$

\begin{prop}\label{propbracket}
Let $X(iu_1,T_1,iv_1)$ and $X(iu_2,T_2,iv_2)$ be two vector fields
in ${\mathfrak g}$. Then the Lie bracket is given by
$$[X(iu_1,T_1,iv_1),X(iu_2,T_2,iv_2)]=X(iu,T,iv),$$ where
\begin{eqnarray*}
u&=&T_1u_2-T_2u_1,\\
T&=&[T_1,T_2]-2( (u_1\Box v_2) -(u_2\Box  v_1)), \\
v&=& T_2^tv_1 - T_1^tv_2\, .
\end{eqnarray*}
\end{prop}

\begin{proof}
The proof is just as is found in \cite[p. 209]{fk}.
\end{proof}
We
note that $\mathfrak{l}$, the Lie algebra of $L$ is
given by
$$\mathfrak{l}=\set{X(0,T,0)\mid \mathfrak{h}\ni T=-T^t},$$
where ${}^t$ denotes the transpose of $T$.

\begin{prop}
The Lie algebra of $K$ is given by
$${\mathfrak k}=\set{ X(iu,T,iu)\mid u\in V, \;T\in {\mathfrak l}\;}\, .$$
\end{prop}

\begin{proof} The map $s(X)=  -iXi$
takes vector fields acting in the upper half plane to those acting
on the right half plane and vice versa.  The vector fields of the
form $X(-u,T,u)$, $T\in \mathfrak{l}$, with the obvious notation,
form the Lie algebra for the group acting on the upper
half plane that fixes $ie$.  (c.f. \cite[p. 210]{fk}). Furthermore,

\begin{eqnarray*}
s(X(-u,T,u))&=&-iX(-u,T,u)i(z)\\
&=&-i(-u+T(iz)-P(iz)u)\\
&=&iu +Tz -P(z)iu\\
&=&X(iu,T,iu)
\end{eqnarray*} and this implies the proposition.

\end{proof}

\begin{prop} The Killing form, $B$, on ${\mathfrak g}$ is given by
  $$B(X(iu_1,T_1,iv_1),X(iu_2,T_2,iv_2))= B_\circ(T_1,T_2) + 2\Tr(T_1T_2) -4\frac{n}{r}(\rip{u_1}{v_2}+\rip{u_2}{v_1}),$$
  where $B_\circ(\cdot,\cdot)$ is the Killing form on ${\mathfrak
  h}$. It is nondegenerate on ${\mathfrak g}$.

\end{prop}
\begin{proof}
c.f.   \cite[p. 28]{s}.
\end{proof}

We now define
$${\mathfrak p}=\set{X(iu,T,-iu)\mid u\in V, \; T=T^t\in
{\mathfrak h}}\, .$$
It is not difficult to see that the Killing form
is negative definite on ${\mathfrak k}$ and positive definite on
${\mathfrak p}$. Moreover,
$${\mathfrak g}={\mathfrak k}\oplus {\mathfrak p}$$
is the Cartan decomposition of $\mathfrak{g}$
corresponding to the Cartan involution
$\Theta:{\mathfrak g}\rightarrow {\mathfrak g}$
given by
$$\Theta(X(iu,T,iv))=X(iv,-T^t,iu)\, .$$

Let ${\mathfrak g}_{\mathbb C}$ be the complexification
of ${\mathfrak g}$ which we
will identify with the set of all vector fields of the form
$X(z,T,w)$, where $z,w \in V_{\mathbb C}$ and
$T\in {\mathfrak h}_{\mathbb C}$. We will let $[\cdot,\cdot]$ denote the complex
linear extension of the bracket given in Proposition
\ref{propbracket}. Specifically, we have
$$[X(z_1,T_1,w_1),X(z_2,T_2,w_2)]=X(z,T,w)\, ,$$
where
\begin{eqnarray*}
z&=&T_1z_2-T_2z_1\\
T&=& [T_1,T_2]+2((z_1\Box w_2) -(z_2\Box w_1))\\
w&=& T_2^tw_1-T_1^tw_2.
\end{eqnarray*}
Let $T\in {\mathfrak h}_{\mathbb C}$ and write $T=T_1+iT_2$ where
$T_1, T_2 \in {\mathfrak h}$.  If $T=T^t$ then $T_1$ and $T_2$ are
likewise self adjoint.  Any self adjoint operator in ${\mathfrak h}$
is a left translation operator $L(x)$, for some $x\in V$.  It
follows then that $T=L(x)+iL(y)=L(x+iy)$, for some $x,y \in V$.
Therefore, the self adjoint operators in ${\mathfrak h}_{\mathbb C}$
are left multiplication operators on $V_{{\mathbb C}}$ by elements
in $V_{\mathbb C}$.

Let $\textsc{z}=X(-e,0,-e)$. Then an easy calculation shows that
$\textsc{z}$ is in the center of ${\mathfrak k}_{\mathbb C}$.
Furthermore, $\ad(\textsc{z})$ has eigenvalues $\pm 2$ on
${\mathfrak p}_{\mathbb C}$, the complexification of ${\mathfrak p}$
in ${\mathfrak g}_{\mathbb C}.$ Indeed, let
\begin{equation}\label{ppm}
{\mathfrak p}_{+} = \set{X(z,L(2z),-z)\mid z\in V_{\mathbb C}} \quad
\text{and} \quad {\mathfrak p}_-=\set{X(z,-L(2z),-z)\mid z\in
V_{\mathbb C}}\, .
\end{equation} Then for $X(z,L(2z),-z)\in {\mathfrak p}_+$
$$[X(-e,0,-e),X(z,L(2z),-z)]=2X(z,L(2z),-z)$$
and for
$X(z,-L(2z),-z)\in {\mathfrak p}_-$
$$[X(-e,0,-e),X(z,-L(2z),-z)]=-2X(z,-L(2z),-z)\, .$$
Since $${\mathfrak p}_{\mathbb C}= {\mathfrak p}_+ \oplus {\mathfrak
p}_-$$ it follows that ${\mathfrak p}_+$ is the $+2$-eigenspace of
$\ad(\textsc{z})$ and ${\mathfrak p}_-$ is the $-2$-eigenspace of
$\ad(\textsc{z})$. Note, that both $\mathfrak{p}_+$ and
$\mathfrak{p}_-$ are Abelian subalgebras of
$\mathfrak{p}_\mathbb{C}$.

\subsection{$L$-fixed vectors}  The group $K$ (and its Lie algebra
${\mathfrak k}$) naturally acts on ${\mathfrak g}_{\mathbb C}.$ We
are interested in the set of vectors, ${\mathfrak g}_{\mathbb C}^L$,
fixed by the action of the subgroup $L$ or, equivalently, those
vectors annihilated by ${\mathfrak l}$ via the adjoint
representation. First of all, since
$${\mathfrak g}_{\mathbb C}= \mathfrak{p}_+\oplus{\mathfrak
k}_{\mathbb C} \oplus {\mathfrak p}_-$$ is a  decomposition into
${\mathfrak k}_{\mathbb C}$-invariant subspaces it follows that
$${\mathfrak g}_{\mathbb C}^L = \mathfrak{p}_+^L \oplus{\mathfrak k}_{\mathbb C}^L
 \oplus{\mathfrak p}_-^L.$$  Let
$\textsc{x}=\frac{1}{2} X(e,2L(e),-e)$, $\textsc{y}=\frac{1}{2}
X(-e,2L(e),e)$, and $\textsc{z}=X(-e,0,-e)$ as above. Then
$\textsc{x}\in \mathfrak{p}_+$, $\textsc{y}\in {\mathfrak p}_-$, and
$\textsc{z} \in {\mathfrak k}_{\mathbb C}$ and each are fixed by
$L$. Furthermore, if $\mathfrak{s}$ is the ${\mathbb C}$-span of
$\set{\textsc{x},\textsc{y},\textsc{z}}$ then ${\mathfrak s}$ is a
Lie subalgebra isomorphic to ${\mathfrak sl}(2,{\mathbb C})$.
Indeed, we need only  observe that

\begin{eqnarray*}
\left[\textsc{x},\textsc{y}\right]&=& \textsc{z},\\
\left[\textsc{z}, \textsc{x}\right]&=& 2\textsc{x},\\
\left[\textsc{z},\textsc{y}\right]&=& -2\textsc{y}.
\end{eqnarray*}

\begin{prop} \label{props}
With the notation established above we have
$${\mathfrak g}_{\mathbb C}^L={\mathfrak s}\, .$$
\end{prop}

\begin{proof}
It follows by \cite{ho}, Theorem 1.3.11, that $\dim_{\mathbb{R}} \mathfrak{g}^L=
\dim_{\mathbb{C}}\mathfrak{g}_\mathbb{C}^L=3$. The claim
follows as
$\dim_{\mathbb{C}}\mathfrak{s}=3$.
\end{proof}

\section{Highest weight representations of $G$}\label{s3}
\noindent In this section we will discuss a well know series of
representations of $G$ acting on  spaces,
$\mathcal{H}_\nu(T(\Omega))$, of holomorphic functions defined on
the tube domain $T(\Omega)= \Omega+iV$.  The Laplace transform,
$\mathcal{L}_\nu$, is a unitary isomorphism of
$L^2(\Omega,d\mu_\nu)$ onto $\mathcal{H}_{\nu}(T(\Omega))$.  We use
this isomorphism to define an equivalent representation of $G$ on
$L^2(\Omega,d\mu_\nu)$.

\subsection{Representation on $\mathcal{H}_{\nu}(T(\Omega))$}
Let $\tilde{G}$ be the universal covering group of $G$
and $\kappa :\tilde{G}\to G$ the covering map.
Then $\tilde{G}$ acts on $T(\Omega)$ via the covering map,
i.e., $g\cdot z=\kappa(g)z$. For
$\nu>1+n(r-1)$ let $\mathcal{H}_{\nu}(T(\Omega))$ be the space of
holomorphic functions $F:T(\Omega)\rightarrow\mathbb{C}$ such that
\begin{equation*}
\left|  \left|  F\right|  \right|
_{\nu}^{2}:=\alpha_{\nu}\int_{T(\Omega
)}|F(ix+y)|^{2}\Delta(y)^{\nu-2n/r}\,dxdy<\infty,\label{eq-normub}%
\end{equation*}
where
\begin{equation*}
\alpha_{\nu}=\frac{2^{r\nu}}{(4\pi)^{n}\Gamma_{\Omega}(\nu-n/r)}\,.
\end{equation*} (See section \ref{sectiongamma} for the definition of
$\Gamma_\Omega$.)  Then $\mathcal{H}_{\nu}(T(\Omega))$ is a
non-trivial Hilbert space. For $\nu\leq1+n(r-1)$ this space reduces
to $\{0\}$. If $\nu=2n/r$ this is the \textit{Bergman space}. For
$g\in\tilde{G}$ and $z\in T(\Omega)$, let $J(g,z)$ be the
\textit{complex} Jacobian determinant of the action of
$g\in\tilde{G}$ on $T(\Omega)$ at the point $z$.  Then
\[
J(ab,z)=J(a,b\cdot z)J(b,z)
\]
for all $a,b\in\tilde{G}$ and $z\in T(\Omega)$.  Recall that
the \textit{genus} of $T(\Omega )$ is
$p=\frac{2n}{r}$.  It is well known that for $\nu>1+n(r-1)$ that
\begin{equation*}
\pi_{\nu}(g)F(z)=J(g^{-1},z)^{\nu/p}F(g^{-1}\cdot z)
\end{equation*}
defines a unitary irreducible representation of $\tilde{G} $, c.f.
\cite{rv}, \cite{w} and \cite{fk}, for example.

For the following we need the explicit form of $J(g,z)$ on
the generators $\tau_{iu}$, $\rho_T$, and $\sigma_{iv}$.

\begin{lem}\label{le-Jgen} The multiplier, $J$, satisfies
$$ \begin{array}{lll}
  (1) & J(\tau_{iu},z)=1,  & u\in V;\\[.5ex]
  (2) & J(\rho_T,z)=\Det T, & T\in {\mathfrak h};\\[.5ex]
  (3) & J(\sigma_{iv},z)=\det (e+izv)^{-p}, & v\in V.\\[1ex]
\end{array}$$
\end{lem}

\begin{proof} In the following $w$ denotes arbitrary element
in $V_\mathbb{C}$ and $z\in T(\Omega )$.
\smallskip

\noindent
(1) Let $u\in V$. Then
$$\frac{d}{dt}\tau_{iu}(z+tw)|_{t=0}=\frac{d}{dt}(z +tw +iu)|_{t=0}=w.$$
Hence $J(\tau_{iu},z)=1$.
\smallskip

\noindent
(2) Let $T\in H$. We then have
$\frac{d}{dt}\rho_T(z+tw)|_{t=0}=Tw$.
Hence $J(\rho_T,z)=\Det T.$
\smallskip

\noindent
(3) For $v\in V$ we get by Lemma \ref{le-DI} and the chain rule
$$\frac{d}{dt}\sigma_{iv}(z+tw)|_{t=0}=\frac{d}{dt}((z+tw)^{-1}+iv)^{-1}|_{t=0}=
  (P(z)P(z^{-1}+iv))^{-1}w\, .$$
Hence $J(\sigma_{iv},z)=((\det z)\det(z^{-1}+iv))^{ -2n/r}
=\det(e+izv)^{-2n/r }$.
%\end{enumerate}
\end{proof}

Recall that the space of smooth vectors $\mathcal{H}_{\nu}(T(\Omega))^\infty$
in $\mathcal{H}_{\nu}(T(\Omega))$ is the space of all
$F\in \mathcal{H}_{\nu}(T(\Omega))$ such
that
$$\mathbb{R}\ni t\mapsto \pi_\nu (\exp tX)F\in \mathcal{H}_{\nu}(T(\Omega))$$
is smooth for all $X\in \mathfrak{g}$.
We denote also by $\pi_{\nu}$  the action of the Lie algebra ${\mathfrak
g}$ and the  complex linear extension to ${\mathfrak g}_{\mathbb
C}$. For $F\in \mathcal{H}_{\nu}(T(\Omega))$ and $X\in \mathfrak{g}$
this action is given by
$$\pi_\nu (X)F=\lim_{t\to 0}\frac{\pi_\nu(\exp tX)F-F}{t}\, .$$
If $F$ is a complex valued holomorphic function on $T(\Omega)$ we
let, as before,  $$D_wF(z)=DF(z)w=\frac{d}{dt}F(z+tw)|_{t=0}$$ be
the (non normalized) directional derivative of $F$ in the direction
$w\in V_{\mathbb C} $. As the point evaluation maps $F\mapsto F(z)$
are continuous linear functionals in $ \mathcal{H}_{\nu}(T(\Omega))$
it follows easily that
\begin{eqnarray*}
\pi_\nu (X)F(z)&=&  \frac{d}{dt} J(\exp (- t X),z)^{ \nu /p }F(\exp
(- t X)z)|_{t=0}
\\
&=& J(\exp (-t X),z)^{\nu/p }|_{t=0}F(z) +\frac{d}{dt}F( \exp (-t X
) z)|_{t=0},
\end{eqnarray*}
for all $z\in T(\Omega )$, $X\in \fg_{\mathbb{C}}$, and $F\in
\mathcal{H}_{\nu}(T(\Omega))^\infty$. The following proposition
gives the action of $\mathfrak {n}^+$, $\mathfrak {n}^-$ and
$\mathfrak{h}$, and hence the full Lie algebra,
on $\mathcal{H}_{\nu}(T(\Omega))$:

\begin{prop} Let $F\in \mathcal{H}_{\nu}(T(\Omega))^\infty$. Then
the subalgebras ${\mathfrak n^+}$, ${\mathfrak h}$, and ${\mathfrak
n^-}$ act by the following formulas:
$$\begin{array}{lll}
(1) & \pi_{\nu} (X(iu,0,0))F(z) = -D_{iu}F(z),&
X(iu,0,0)\in\mathfrak{n}^+;\\[.5ex]
(2) & \pi_{\nu}(X(0,T,0))F(z)=-\frac{\nu}{p}\tr(T)F(z)-D_{Tz}F(z),&
 X(0,T,0)\in\mathfrak{h};\\[.5ex]
(3)& \pi_{\nu}(X(0,0,iv))F(z)=i\nu\tr(zv)F(z)+D_{P(z)iv}F(z),&
X(0,0,iv)\in\mathfrak{n}^-.\\[1ex]
\end{array}$$
\end{prop}

\begin{proof}
(1) Let $u\in V$. The formula $t\mapsto \tau_{tiu}(z)=z+itu$ defines
the one parameter subgroup in the direction $X(iu,0,0)\in {\mathfrak
n}^+$. By Lemma \ref{le-Jgen}, we then have
\begin{eqnarray*}
\pi_v(X(iu,0,0))F(z) &=& \frac{d}{dt}
J(\tau_{tiu}^{-1},z)^{\nu /p }F(\tau_{tiu}^{-1}z)|_{t=0}\\
&=&\frac{d}{dt}J(\tau_{-tiu},z)^{\nu /p }|_{t=0} F(z) +
\frac{d}{dt}F( \tau_{-tiu} z)|_{t=0}\\
&=&D_{-iu}F(z).
\end{eqnarray*}
\smallskip

\noindent
(2) Let $T\in {\mathfrak h}$. Then $t\mapsto \rho_{\exp
  (tT)}(z)=\exp(tT) z$ defines the one parameter subgroup in the
direction $X(0,T,0)\in {\mathfrak h}$. By Lemma \ref{le-Jgen}, we
have
\begin{eqnarray*}
\pi_\nu(X(0,T,0))F(z)&=& \frac{d}{dt}J(\rho_{\exp(-tT)},z)^{\nu /p }
F(\rho_{\exp(-tT)}z)|_{t=0}\\
&=& \frac{d}{dt}J(\rho_{\exp(-tT)},z)^{\nu /p}|_{t=0} F(z) +
\frac{d}{dt}F(\rho_{\exp(-tT)}z)|_{t=0}\\
&=& \frac{-\nu}{p}\tr(T)F(z) -D_{Tz}F(z).
\end{eqnarray*}
\smallskip

\noindent (2) Finally, let $v\in V$. Then $z\mapsto \sigma_{tiv}(z)=
(z^{-1}+itv)^{-1}$ defines the one parameter subgroup  in the
direction $X(0,0,iv)$. Again, by Lemma \ref{le-Jgen} and
by Lemma \ref{le-DI}, we  have
\begin{eqnarray*}
\pi_\nu(X(0,0,iv))F(z)&=& \frac{d}{dt}
J(\sigma_{itv}^{-1},z)^{\nu /p }F(\sigma_{itv}^{-1}z)|_{t=0}\\
&=& \frac{d}{dt}J(\sigma_{-itv},z)^{\nu /p}|_{t=0}F(z) +
\frac{d}{dt}F(z^{-1}-i t v)^{-1})|_{t=0}\\
&=&\frac{d}{dt} \det(e-i t z v)^{-\nu}|_{t=0}F(z) +
D F(z)[-P(z^{-1})^{-1}(-iv)]\\
&=& i\nu\tr(z v)F(z)+D_{P(z)iv}F(z).\\
  \end{eqnarray*}
\end{proof}

\subsection{The Laplace Transform}
Let $L^2(\Omega,d\mu_\nu)$ be the Hilbert space of all Lebesgue
measurable functions on $\Omega$ such that
$$\norm{f}^2=\int_{\Omega}^{} \abs{f(x)}^2\, d\mu_\nu(x)<\infty,$$
where $d\mu_\nu(x)=\Delta^{\nu-n/r}(x)\; dx$.  For $f\in
L^2(\Omega,d\mu_\nu)$  the Laplace transform is defined by the formula
$$\mathcal{L}_\nu(f)(z)=\int_{\Omega}^{ }e^{-\rip{z}{x}}f(x)\,
d\mu_\nu(x)\, .$$

\begin{prop} \label{ltunitary} Let $f\in L^2(\Omega,d\mu_\nu)$.
Then $\mathcal{L}_\nu f
\in \mathcal{H}_{\nu}(T(\Omega))$.  Furthermore, the map
$$\mathcal{L}_\nu: L^2(\Omega,d\mu_\nu)\rightarrow
\mathcal{H}_{\nu}(T(\Omega))$$ is a unitary isomorphism.
\end{prop}
\begin{proof}
C.f. \cite{rv, doz2}.
\end{proof}

\subsection{Representation on  $L^2(\Omega,d\mu_\nu)$}
By Proposition \ref{ltunitary} we can define an equivalent
representation, $\lambda_\nu$,  of $G$ on $L^2(\Omega,d\mu_\nu)$  so
that $\mathcal{L}_\nu$ is an intertwining operator. Specifically,
$$\lambda_\nu(g)f=\mathcal{L}_\nu^{-1} \pi_\nu(g) \mathcal{L}_\nu f,$$ for $g\in
G$. We denote by $L^2(\Omega,d\mu_\nu)^\infty$ the space of smooth
vectors. As usual we will let $\lambda_\nu$ also denote the action
of the Lie algebras  ${\mathfrak g}$ and $\mathfrak{g}_{\mathbb C} $
on $L^2(\Omega,d\mu_\nu)^\infty$. Note that this representation is
not geometric in the sense that  $\tilde G$ does  not act naturally
on $\Omega$,  only the subgroup - with the obvious notation -
$\widetilde{G(\Omega)}\cap G$ acts on $\Omega$. We follow \cite{fk}
to define the Bessel operator $\mathcal{B}_\nu:
C^{\infty}(V)\rightarrow C^{\infty}(V) \otimes V_{\mathbb C}$
formally by
\begin{equation}\label{eq-Bessel}
\mathcal{B }_\nu=P\left(\frac{\partial}{\partial x}\right)x + \nu
\frac{\partial}{\partial x}\, .
\end{equation}
If $\set{e_i}_{i=1}^{ n}$ is an orthonormal basis of $V$ and
$(x_1,\ldots ,x_n)$ the corresponding coordinate functions, then
$$\mathcal{B}_\nu f(x) =\sum_{i,j}^{ }
\frac{\partial^2 f }{\partial x_i \partial x_j}(x)P(e_i,e_j)x + \nu
\sum_{i} \frac{\partial f}{\partial x_i}(x)e_i\, .$$ By the
definition given in Equation \ref{eq-Bessel} this formula is
obviously basis independent. We refer to  \cite[p. 322]{fk} for more
details. For $w\in V_\mathbb{C}$ we define the differential operator
$\mathcal{B}_{\nu,w}$ by
\begin{equation}\label{eq-Bnuw}
\mathcal{B}_{\nu ,w}f(x)=\rip{\mathcal{B}_\nu f(x)}{w}\, .
\end{equation}

\begin{thm}\label{th-3.4}
Let $f\in L^2(\Omega,d\mu_\nu)^\infty$. The representation $\lambda_\nu$
of ${\mathfrak g}$ is described as follows:
$$\begin{array}{lll}
(1)& \lambda_\nu(X(iu,0,0))f(x)= \tr(iux)f(x)\, , & X(iu,0,0)\in
{\mathfrak n}^+;\\[.5ex]
(2) & \lambda_\nu(X(0,T,0))f(x)= \frac{\nu}{p}\Tr(T)f(x)+
D_{T^tx}f(x)\, , &  X(0,T,0)\in {\mathfrak h};\\[.5ex]
(3) & \lambda_\nu(X(0,0,iv))f(x) = -\mathcal{B}_{\nu,iv} f(x)\, , &
X(0,0,iv)\in {\mathfrak  n}^-.\\[1ex]
\end{array}$$
\end{thm}
\begin{proof}
%\begin{enumerate}
(1) Let $u\in V$ and for convenience let $m=\nu-n/r$. Let
$w=\frac{1}{2}z\in \Omega$. Then $w+\Omega$ is an open neighborhood
of $z$ and for $f\in L^2(\Omega,d\mu_\nu)^\infty$ we have
$$| e^{-\rip{y}{x}}\rip{iu}{x}f(x)\Delta^m(x)|\le
|e^{-\rip{w}{x}}\rip{iu}{x}f(x)\Delta^m(x)|,$$ for all $y\in
w+\Omega$. As $e^{-\rip{w}{\cdot }}\rip{iu}{\cdot }f\Delta^m\in
L^1(\Omega ,d\mu_\nu)$ we can interchange the integration and
differentiation in the following to get
\begin{eqnarray*}
\pi_\nu(X(iu,0,0))\mathcal{L}_\nu f(z)&=&
-D_{iu}\mathcal{L}_\nu f(z)\\
 &=&
-\int_{\Omega}^{ }D_{iu} ( e^{-\rip{z}{x}} )f(x)\Delta^m(x) \,
dx\\
&=& \int_{\Omega}^{ }e^{-\rip{z}{x}}\rip{iu}{x}f(x)\Delta^m(x)\,
dx\\
&=&\mathcal{L}_\nu(\tr(iux)f(x))(z).
\end{eqnarray*}
\smallskip

\noindent
(2)
In \cite[p. 191]{doz2} we determine the action of $H$ on
$L^2(\Omega,d\mu_\nu)$ as follows:
$$\lambda_\nu(h)f(x)=\Det(h)^{\frac{\nu}{p}}f(h^t x),$$  $h\in
H$.  Differentiation of this formula gives (2).
\smallskip

\noindent
(3) By Proposition XV.2.4 of \cite{fk} we have
\begin{eqnarray*}
\mathcal{L}_\nu(\mathcal{B}_{\nu,iv} f)(z)&=&
-\rip{P(z)\frac{\partial}{\partial z}+ \nu z}{iv}\mathcal{L}_\nu
f(z)\\
&=&-\left( D_{P(z)iv}+\nu \tr(izv)\right)\mathcal{L}_\nu f(z)\\
&=& -\pi_\nu(X(0,0,iv))\mathcal{L}_\nu f(z),
\end{eqnarray*} from which the result follows.
%\end{enumerate}
\end{proof}

\begin{rem}\label{remaction}
Each of these formulas extend to the complexification in an obvious
way:
$$\begin{array}{lll}
(1)&\lambda_\nu(X(w,0,0))f(x)=\tr(wx)f(x), &  X(w,0,0)\in {\mathfrak
n^+_{\mathbb C} };\\[.5ex]
(2)&\lambda_\nu(X(0,T,0))f(x)= \frac{\nu}{p}\Tr(T)f(x)+
  D_{T^tx}f(x), &  X(0,T,0)\in {\mathfrak h}_{\mathbb C};\\[.5ex]
(3)& \lambda_\nu(X(0,0,w))f(x) = - \mathcal{B}_{\nu ,w} f(x),&
X(0,0,w)\in {\mathfrak n}^-_{\mathbb C}.\\[1ex]
\end{array}$$
\end{rem}

\section{Laguerre functions}\label{s4}
\noindent We continue with the assumption that $V$ is a simple
Euclidean Jordan algebra with rank $r$, dimension $n$, and degree
$d$; c.f. \cite[p. 71 and 98]{fk} for the definition of the degree
of a Jordan algebra.  Let $c_i$, $i=1,\ldots, r$, be a fixed Jordan
frame and $V^{(k)}$ the $+1$ eigenspace of the operator $L(c_1 +
\cdots +c_k)$.  Then each $V^{(i)}$ is a Jordan subalgebra of $V$
and we have the following inclusions:
$$V^{(1)}\subset V^{(2)}\subset \cdots \subset V^{(r)}=V.$$
Let $\det_i$, $i=1,\ldots, r$ be the determinant function on
$V^{(i)}$ and  define $\Delta_i(x) = \det(P_i x)$, where $P_i$ is
orthogonal projection of $V$ onto $V^{(i)}.$ These are the principal
minors, they are homogenous polynomials of degree $i$, and
$\Delta_r(x)=\det x$. For convenience we write $\Delta=\Delta_r$,
c.f. \cite[p. 114]{fk} for details. For ${\mathbf
s}=(s_1,\ldots,s_n)\in {\mathbb C}^r$ define
$$\Delta_{\mathbf s}=\Delta_1^{s_1-s_2}\Delta_2^{s_2-s_3}\cdots
\Delta_r^{s_r}.$$
 For ${\mathbf
m}=(m_1,\ldots,m_r)$ a sequence on nonnegative integers we write $
{\mathbf m}\ge 0$ to mean $m_1\ge m_2 \ge \cdots \ge m_r\ge 0$. Let
$$\Lambda=\set{{\mathbf m}\mid {\mathbf m}\ge 0}.$$ Then $\Delta_{\mathbf m}$
are the \textit{generalized power functions} of degree
$\abs{\mathbf{m} }=m_1 + \cdots + m_r$. Since $\Delta_{\mathbf m}$
is a polynomial function on $V$ it extends uniquely to a holomorphic
polynomial function on $V_{\mathbb C}$.

\subsection{$L$-invariant Polynomials}  We define $\psi_{\mathbf m}$
by the following formula $$\psi_{\mathbf m}(x)=\int_{L}^{
}\Delta_{\mathbf m}(lx)\, dl,\quad x\in V,$$  where $dl$ is
normalized Haar measure on $L$. The function $\psi_{\mathbf m}$ is a
nonzero  $L$-invariant polynomial on $V$, for each ${\mathbf m}\in
\Lambda$, which also extends uniquely to a holomorphic function on
$V_{\mathbb C}$. Furthermore, the set of $L$-invariant polynomials
is spanned by the set of all $\psi_{\mathbf m}, \; {\mathbf m}\in
\Lambda$. Moreover, if $\mathcal{P}_k(V)$ denotes the set of
$L$-invariant polynomials on $V$ of degree at most $k$ then
$\mathcal{P}_k(V)$ is spanned by all those $\psi_{\mathbf m}$ with
$\abs{\mathbf{m} }\le k$. The function $\psi_{\mathbf m}(e+x)$ is
also an $L$-invariant polynomial of degree $\abs{\mathbf{m} }$ and
has an expansion that defines the generalized binomial coefficients
$\binom{{\mathbf m} }{{\mathbf n}}$:
$$\psi_{\mathbf m}(e+x)=\sum_{\abs{{\mathbf n}}\le \abs{\mathbf m} }^{ }
\binom{{\mathbf m}}{{\mathbf n}} \psi_{\mathbf n}(x)\, .$$

\subsection{The Gamma Function} \label{sectiongamma} For convenience we
also reproduce the gamma function of
the symmetric cone $\Omega$: Let ${\mathbf s}\in {\mathbb C}^r$ and
define $$\Gamma_\Omega({\mathbf s})=\int_{\Omega}^{ }e^{-\tr
x}\Delta_{\mathbf s}(x) \Delta^{-n/r}(x)\; dx,$$
where
$\Delta(x)=\Delta_r(x)$ as before. For $\nu$ a real number  and ${\mathbf
m}\in \Lambda$ define
$$(\nu)_{\mathbf m}=\frac{\Gamma_\Omega(\nu + {\mathbf
m})}{\Gamma_\Omega(\nu)},$$ where $\nu+{\mathbf m}$ means to add
$\nu$ to each component of ${\mathbf m}$.

\subsection{The Generalized Laguerre Functions} The Laguerre
polynomials are defined by
\begin{equation}\label{eq-Lpol}
L^\nu_{\mathbf m}(x)=(\nu)_{\mathbf m}
\sum_{ \abs{{\mathbf n}}\le \abs{\mathbf m}}^{ } \binom{{\mathbf
m}}{{\mathbf n}} \frac{1}{(\nu)_{\mathbf n} }\psi_{\mathbf n}(-x)\, ,
\end{equation}
and the generalized Laguerre functions by
\begin{equation}\label{eq-Lfct}
\ell_{\mathbf m}^\nu (x)=e^{-\tr x}L^\nu_{\mathbf m}(2x)\, ,
\end{equation} c.f. \cite{fk}, p. 343.

\begin{rem} Let ${\mathbf 0}$ denote the multiindex with entries $0$.
Then $\ell_{\mathbf 0}^\nu=e^{-\tr x}$.  Let $X(z,2L(z),-z)\in
{\mathfrak p}_+$, c.f. Equation \ref{ppm}. Then a straightforward
calculation gives $$\lambda_\nu(X(z,2L(z),-z))\ell_{\mathbf
0}^\nu=0,$$ for all $z\in V_{\mathbb C}$.  Thus $\ell_{\mathbf
0}^\nu$ is the highest weight vector for $\lambda_\nu$.
\end{rem}

\begin{thm}
Let $L^2(\Omega,d\mu_\nu)^L$ be the space of $L$-invariant function
in $L^2(\Omega,d\mu_\nu)$. Then the Laguerre functions form an
orthogonal basis of $L^2(\Omega,d\mu_\nu)^L$. Moreover,
$$\norm{\ell_{\mathbf m}^\nu (x)}^2 = \frac{1}{2^{r\nu}}\frac{1}{d_{\mathbf m}
} \left(\frac{n}{r}\right)_{\mathbf m} \Gamma_\Omega(\nu+{\mathbf
m}).$$
\end{thm}
\begin{proof}
This is Theorem 4.1 in \cite{doz2}. See also \cite[p. 344]{fk}.
\end{proof}

\section{Differential Recursion Relations}\label{s5}
\noindent Recall that $ {\mathfrak g}_{\mathbb C} ^L$ is the set of
vector fields in ${\mathfrak g}_{\mathbb C} $ invariant under the
adjoint action of $L$. Proposition \ref{props} establishes that
${\mathfrak s}={\mathfrak g}_{\mathbb C} ^L$ is isomorphic to
${\mathfrak sl}(2,{\mathbb C})$ and is spanned be the vector fields
\begin{eqnarray*}
\textsc{x}&=&\frac{1}{2} X(e,2L(e),-e) \in {\mathfrak p}^+\, ,\\
\textsc{y}&=&\frac{1}{2} X(-e,2L(e),e) \in {\mathfrak p}^-\, ,\\
\textsc{z}&=&X(-e,0,-e) \in {\mathfrak k}_{\mathbb C}\, .
\end{eqnarray*}
Our main theorem generalizes the classical differential recursion
relations on Laguerre functions by way of the explicit action of
${\mathfrak s}$ on $L^2(\Omega,d\mu_\nu)^{\infty L}$. Set $B_\nu =
\mathcal{B}_{\nu, e}$.
\begin{prop}\label{propsaction} Let $f\in L^2(\Omega,d\mu_\nu)^\infty$.
With notation as above we have
\begin{enumerate}
\item $\lambda_\nu(\textsc{x})f(x)=\frac{1}{2}(\tr x + r\nu + 2
D_x + B_{\nu})f(x)\, ,$
\item $\lambda_\nu(\textsc{y})f(x)=\frac{1}{2}(-\tr x + r\nu + 2
D_x - B_{\nu})f(x)\, ,$
\item  $\lambda_\nu(\textsc{z})f(x)= (-\tr x + B_{\nu})f(x)\, .$
\end{enumerate}
\end{prop}

\begin{proof}
These formulas follow directly from Remark \ref{remaction}.
\end{proof}

\begin{thm}\label{th-5.2}
The Laguerre functions are related by the following differential
recursion relations:
\begin{enumerate}
\item $(-\tr x + B_\nu )\ell_{\mathbf m}^\nu
  (x)=-(r\nu +2 \abs{{\mathbf m}})\ell_{\mathbf m}^\nu (x)$,
\item $(\tr x +r\nu +2D_x +B_\nu )\ell_{\mathbf m}^\nu(x)=-2\sum_{j=1
}^{r}\binom{\mathbf m}{{\mathbf
m}-\gamma_j}(m_j-1+\nu-(j-1)\frac{d}{2})\ell_{{\mathbf
m}-\gamma_j}^\nu(x)$,
\item $(\tr x -r\nu -2D_x +B_\nu )\ell_{\mathbf
m}^\nu(x)=-2\sum_{j=1}^{r }c_{\mathbf m}(j)\ell_{{\mathbf
m}+\gamma_j}^\nu(x)$,
\end{enumerate}
where the constants $c_{\mathbf m}(j)$  are defined by
$$c_{\mathbf m}(j)=\prod_{k\ne
j}\frac{m_k-m_j-\frac{d}{2}(k-j+1)}{m_k-m_j-\frac{d}{2}(k-j)}\, .$$
\end{thm}

\begin{rem}
In the formulas above $\mathbf{ \gamma}_j$ is the multiindex with
$1$ in the $j^{\text{th}}$ position and $0$'s elsewhere. It should
be understood that if $\mathbf{m + \gamma_j}$ or
$\mathbf{m-\gamma_j}$ is not in $\Lambda$ then the corresponding
Laguerre function does not appear in the sum.
\end{rem}

\begin{proof}
Let $\xi=X(e,0,e)=e-P(z)e$. Then $\xi=-\textsc{z}$ is the vector
field given by the same symbol in \cite{doz2}.  By Lemma 5.5 of
\cite{doz2} $$\pi_\nu(\xi)q_{{\mathbf m},\nu}=(r\nu+2\abs{{\mathbf
m}})q_{{\mathbf m},\nu},$$ where $\Gamma_\Omega({\mathbf m}+
\nu)q_{{\mathbf m},\nu}= \mathcal{L}_{\nu}(\ell_{\mathbf m}^\nu)$.
By the unitary equivalence of $\pi_\nu$ and $\lambda_\nu$ we
correspondingly have
$$\lambda_\nu(\xi)\ell_{\mathbf m}^\nu=(r\nu+2\abs{{\mathbf m}})\ell_{\mathbf
m}^\nu.$$ On the other hand,
$$\lambda_\nu(\xi)\ell_{\mathbf m}^\nu=-\lambda_\nu(\textsc{z})\ell_{\mathbf m}^\nu=(\tr x
- B_\nu)\ell_{\mathbf m}^\nu$$ by Proposition \ref{propsaction}.
Part (1) now follows.

Let
$$L_k^2(\Omega,d\mu_\nu)=\set{f\in L^2(\Omega,d\mu_\nu)^\infty \mid
\lambda_\nu(z)f =-(r\nu+2k)f}\, .$$
Since $\lambda_\nu$ is an
irreducible highest weight representation it is well known that
$L_k^2(\Omega,d\mu_\nu)$ is finite dimensional, nonzero if $k\ge 0$,
and $$L^2(\Omega,d\mu_\nu)=\bigoplus_{k=0}^\infty L_k^2(\Omega,d\mu_\nu)\, .$$
Furthermore, part (1) implies that $\ell_{\mathbf m}^\nu\in
L_{\abs{{\mathbf m}}}^2(\Omega,d\mu_\nu).$  For
$\textsc{w}=X(w,2L(w),-w)\in {\mathfrak p}^+$ and $f\in
L_k^2(\Omega,d\mu_\nu)$ we have
\begin{eqnarray*}
\lambda_\nu(\textsc{z})\lambda_\nu(\textsc{w})f&=&\lambda_\nu(\textsc{w})\lambda(\textsc{z})f
+\lambda_\nu([\textsc{z},\textsc{w}])f\\
&=& -(r\nu+2k)\lambda_\nu(\textsc{w})f + 2\lambda_\nu(\textsc{w})f\\
&=& -(r\nu + 2(k-1))\lambda_\nu(\textsc{w})f.
\end{eqnarray*} This implies that $\lambda_\nu(\textsc{w})f\in
L_{k-1}^2(\Omega,d\mu_\nu)$. Similarly, for $\textsc{w}\in
{\mathfrak p}^-$, we have $\lambda_\nu(\textsc{w})f\in
L_{k+1}^2(\Omega,d\mu_\nu)$.

Now let $Z_0=\frac{1}{2}(\textsc{x}+\textsc{y})$. Then $Z_0=X(0,I,0)$
is the Euler vector field $z\frac{\partial}{\partial z}$ given in
\cite[p.161]{doz2}. By Theorem 7.9 of \cite{doz2} (and its proof) we
have $$-2\lambda_\nu(Z_0)\ell_{\mathbf m}^\nu=\sum_{j=1
}^{r}\binom{\mathbf m}{{\mathbf
m}-\gamma_j}(m_j-1+\nu-(j-1)\frac{a}{2})\ell_{{\mathbf
m}-\gamma_j}^\nu-\sum_{j=1}^{r }c_{\mathbf m}(j)\ell_{{\mathbf
m}+\gamma_j}^\nu.$$ If $P_k$ denotes orthogonal projection of
$L^2(\Omega,d\mu_\nu)$ onto $L_k^2(\Omega,d\mu_\nu)$ then
$$-\lambda_\nu(\textsc{x})\ell_{\mathbf m}^\nu=P_{\abs{{\mathbf
m}}-1}(-2\lambda_\nu(Z_0)\ell_{\mathbf m}^\nu)=\sum_{j=1
}^{r}\binom{\mathbf m}{{\mathbf
m}-\gamma_j}(m_j-1+\nu-(j-1)\frac{a}{2})\ell_{{\mathbf
m}-\gamma_j}^\nu $$ and
$$-\lambda_\nu(\textsc{y})\ell_{\mathbf m}^\nu=P_{\abs{{\mathbf
m}}+1}(-2\lambda_\nu(Z_0)\ell_{\mathbf m}^\nu)=-\sum_{j=1
}^{r}c_{\mathbf m}(j) \ell_{{\mathbf m}+\gamma_j}^\nu.$$ We obtain
formulas (2) and (3) by again applying Proposition
\ref{propsaction}.

\end{proof}

\begin{rem}
Observe that one-half the difference between formula (2) and (3)
gives Theorem 7.9 of \cite{doz2}.
\end{rem}

\end{document}